\documentclass[10pt]{amsart}
\usepackage{amssymb, amsthm, amsmath}
\usepackage{xy}
\usepackage{epsfig}
\usepackage{psfrag}

\psfrag{lambda}{$\lambda$}
\theoremstyle{plain}
\newtheorem{prop}{Proposition}
\newtheorem{thm}{Theorem}
\newtheorem*{thm1a}{Theorem 1$'$}
\newtheorem*{giamI}{Giambelli I}
\newtheorem*{giamII}{Giambelli II}
\newtheorem{cor}{Corollary}
\newtheorem{lemma}{Lemma}

\theoremstyle{remark}
\theoremstyle{definition}
\newtheorem*{remark}{Remark}

\newtheorem{exa}{Example}
\newcommand{\refthm}[1]{Theorem~\ref{#1}}
\newcommand{\refprop}[1]{Proposition~\ref{#1}}
\newcommand{\reflemma}[1]{Lemma~\ref{#1}}
\newcommand{\refcor}[1]{Corollary~\ref{#1}}

\newcommand{\rank}{\operatorname{rank}}
\newcommand{\tto}{\twoheadrightarrow}
\newcommand{\Schub}{{\mathfrak S}}
\newcommand{\uu}[1]{\underline{\underline{#1}}}
\newcommand{\Gr}{\operatorname{Gr}}
\newcommand{\col}{\operatorname{col}}

\newcommand{\Z}{{\mathbb Z}}

\newcommand{\C}{{\mathbb C}}

\newcommand{\X}{{\mathfrak X}}

\newcommand{\AS}{{\mathfrak S}}

\newcommand{\la}{\lambda}

\newcommand{\al}{\alpha}
\newcommand{\bull}{{\scriptscriptstyle \bullet}}

\DeclareMathOperator{\GL}{GL}

\newcommand{\gequ}{\geqslant}
\newcommand{\lequ}{\leqslant}

\newcommand{\ra}{\rightarrow}

\newcommand{\noin}{\noindent}
\newcommand{\wt}{\widetilde}

\begin{document}

\title[Schubert polynomials and quiver formulas]
{Schubert polynomials and quiver formulas}
\author{Anders S.~Buch, Andrew Kresch, 
Harry Tamvakis, and Alexander Yong}
\date{November 19, 2002}
\thanks{The authors were supported in part by NSF Grant DMS-0070479
  (Buch), an NSF Postdoctoral Research Fellowship (Kresch),
  and NSF Grant DMS-0296023 (Tamvakis).}
\address{Matematisk Institut, Aarhus Universitet, Ny Munkegade, 8000
  {\AA}rhus C, Denmark}
\email{abuch@imf.au.dk}
\address{Department of Mathematics, University of Pennsylvania,
209 South 33rd Street,
Philadelphia, PA 19104-6395, USA}
\email{kresch@math.upenn.edu}
\address{Department of Mathematics, Brandeis University - MS 050,
P. O. Box 9110, Waltham, MA
02454-9110, USA}
\email{harryt@brandeis.edu}
\address{Department of Mathematics,  University of Michigan, 525 East University Ave.,
Ann Arbor, MI 48109-1109, USA}
\email{ayong@umich.edu}

\maketitle 

\section{Introduction}
\label{intro}

The work of Buch and Fulton \cite{BF} established a formula for a
general kind of degeneracy locus associated to an oriented quiver of
type $A$.  The main ingredients in this formula are Schur determinants
and certain integers, the {\em quiver coefficients}, which generalize
the classical Littlewood-Richardson coefficients.  Our aim in this
paper is to prove a positive combinatorial formula for the quiver
coefficients when the rank conditions defining the degeneracy locus
are given by a permutation.  In particular, this gives new expansions
for Fulton's universal Schubert polynomials \cite{F2} and the Schubert
polynomials of Lascoux and Sch\"utzenberger \cite{LS}.

Let $\X$ be a smooth complex algebraic variety and let
\begin{equation} 
\label{eqn:univseq}
G_1 \to \cdots \to G_{n-1} \to G_n \to F_n \to
F_{n-1} \to \dots \to F_1
\end{equation}
be a sequence of vector bundles and morphisms over $\X$,
such that $G_i$ and $F_i$ have rank $i$ for each $i$.  
For every permutation $w$ in the symmetric group $S_{n+1}$ 
there is a degeneracy locus
\[
\Omega_w(G_\bull \to F_\bull)= 
  \{ x \in \X \ | \ \rank(G_q(x) \to F_p(x)) \lequ r_w(p,q) 
  \text{ for all $1\lequ p,q\lequ n$} \},
\]
where $r_{w}(p,q)$ is the number of $i\lequ p$ such that $w(i)\lequ
q$.  The universal double Schubert polynomial $\Schub_w(c;d)$ of
Fulton gives a formula for this locus; this is a polynomial in the
variables $c_i(j)$ and $d_i(j)$ for $1 \lequ i \lequ j\lequ n$. When
the codimension of $\Omega_w(G_\bull \to F_\bull)$ is equal to the
length of $w$, its class $[\Omega_w]$ in the cohomology (or Chow ring)
of $\X$ is obtained by evaluating $\AS_w(c;d)$ at the Chern classes
$c_i(p)$ and $d_i(q)$ of the bundles $F_p$ and $G_q$, respectively.
The quiver formula given in \cite{BF} specializes to a formula for
universal double Schubert polynomials:
\[ \AS_w(c;d) = \sum_\lambda c_{w,\lambda}^{(n)}\, 
   s_{\lambda^1}(d(2)-d(1)) \cdots
   s_{\lambda^n}(c(n)-d(n)) \cdots s_{\lambda^{2n-1}}(c(1)-c(2)) \,.
\]
Here the sum is over sequences of $2n-1$ partitions $\lambda =
(\lambda^1,\dots,\lambda^{2n-1})$ and each $s_{\lambda^i}$ is a Schur
determinant in the difference of the two alphabets in its argument.
The quiver coefficients $c_{w,\lambda}^{(n)}$ can be computed by an
inductive algorithm, and are conjectured to be nonnegative \cite{BF}.

Let $\col(T)$ denote the column word of a semistandard Young tableau
$T$, the word obtained by reading the entries of the columns of the
tableau from bottom to top and left to right.  The following theorem
is our main result. 

\begin{thm}
\label{introthm}
Suppose that $w\in{\mathcal S}_{n+1}$ and
$\lambda=(\lambda^1,\ldots,\lambda^{2n-1})$ is a sequence of
partitions. Then $c_{w,\lambda}^{(n)}$ equals the number of sequences
of semistandard tableaux $(T_{1},\ldots,T_{2n-1})$ such that the shape
of $T_i$ is conjugate to $\lambda^i$, the entries of $T_{i}$ are at
most ${\rm min}(i,2n-i)$, and $\col(T_{1})\cdots\col(T_{2n-1})$ is a
reduced word for $w$.
\end{thm}

As a consequence of Theorem \ref{introthm}, we obtain formulas for the
Schubert polynomials of Lascoux and Sch\"{u}tzenberger, expressing
them as linear combinations of products of Schur polynomials in
disjoint sets of variables. The coefficients in these expansions are
all quiver coefficients $c^{(n)}_{w,\lambda}$; this is true in
particular for the expansion of a Schubert polynomial as a linear
combination of monomials (see, e.g., \cite{BJS,FS}).

The {\em Stanley symmetric functions} or stable Schubert polynomials
\cite{S} play a central role in this paper.  Our main result
generalizes the formula of Fomin and Greene \cite{FG} for the
expansion of these symmetric functions in the Schur basis, as well as
the connection between quiver coefficients and Stanley symmetric
functions obtained in \cite{B}.  In fact, we show that the universal
double Schubert polynomials can be expressed as a multiplicity-free
sum of products of Stanley symmetric functions (Theorem
\ref{thm:stanprd}).

It should be noted that the formula for Schubert polynomials suggested
in this paper is different from the one given in \cite{BF}.  For
example, the formula from \cite[\S2.3]{BF} does not make it clear that the
monomial coefficients of Schubert polynomials are quiver coefficients,
or even that these monomial coefficients are nonnegative.

We review the universal Schubert polynomials and Stanley symmetric
functions in Section \ref{prelims}; in addition, we prove some
required properties.  In Section \ref{qv} we introduce quiver
varieties and we prove Theorem \ref{introthm}.  In the following 
section we apply our results to the case of ordinary Schubert polynomials.
Finally, in Section \ref{giamsec} we use our expressions for single 
Schubert polynomials to obtain explicit Giambelli formulas for the classical
and quantum cohomology rings of partial flag varieties.

The authors thank Sergey Fomin, Peter Magyar, and Alexander Postnikov
for useful discussions.

\section{Preliminaries}
\label{prelims}

\subsection{Universal Schubert polynomials} \label{usp}

We begin by recalling the definition of the double Schubert
polynomials of Lascoux and Sch{\"u}tzenberger \cite{LS, L1}.
Let $X=(x_1,x_2,\ldots)$ and $Y=(y_1,y_2,\ldots)$ be two sequences of
commuting independent variables.  Given a permutation $w \in S_n$, the
double Schubert polynomial $\Schub_w(X;Y)$ is defined recursively 
as follows.  If
$w = w_0$ is the longest permutation in $S_n$ then we set
\[ \Schub_{w_0}(X;Y) = \prod_{i+j \lequ n} (x_i - y_j) \,. \]
Otherwise we can find a simple transposition $s_i = (i,i+1) \in S_n$
such that $\ell(w s_i) = \ell(w) + 1$.  Here $\ell(w)$ denotes the
length of $w$, which is the smallest number $\ell$ for which $w$ can
be written as a product of $\ell$ simple transpositions.  We then
define
\[ \Schub_w(X;Y) = \partial_i(\Schub_{w s_i}(X;Y)) \]
where $\partial_i$ is the divided difference operator given by
\[ \partial_i(f) = 
  \frac{f(x_1,\dots,x_i,x_{i+1},\dots,x_n) - 
        f(x_1,\dots,x_{i+1},x_i,\dots,x_n)}{x_i - x_{i+1}} \,.
\]
The (single) Schubert polynomial is defined by $\Schub_w(X) =
\Schub_w(X;0)$.

Suppose now that $w$ is a permutation in $S_{n+1}$.  If
$u_1,\ldots,u_r$ are permutations, we will write $u_1 \cdots u_r=w$ if
$\ell(u_1) + \cdots + \ell(u_r) = \ell(w)$ and the product of
$u_1,\ldots,u_r$ is equal to $w$. In this case we say $u_1 \cdots u_r$
is a {\em reduced factorization\/} of $w$.

The {\em double universal Schubert polynomial} $\Schub_w(c;d)$ of
\cite{F2} is a polynomial in the variables $c_i(j)$ and $d_i(j)$ for $1
\lequ i \lequ j \lequ n$.  For convenience we set $c_0(j)=d_0(j) = 1$
for all $j$ and $c_i(j) =d_i(j) = 0$ if $i < 0$ or $i > j$.  The
classical Schubert polynomial $\Schub_w(X)$ can be written uniquely in
the form
\[ \Schub_w(X) = \sum a_{i_1,\dots,i_n} \,
   e_{i_1}(x_1)\, e_{i_2}(x_1,x_2) \cdots e_{i_n}(x_1,\dots,x_n) 
\]
where the sum is over all sequences $(i_1,\dots,i_n)$ with each
$i_\alpha \lequ \alpha$ and $\sum i_\alpha = \ell(w)$~\cite{LS3}.  The
coefficients $a_{i_1,\dots,i_n}$ are uniquely determined
integers depending on $w$.
Define the single universal Schubert polynomial for $w$ by
\[ \Schub_w(c) = \sum a_{i_1,\dots,i_n} \, c_{i_1}(1)\, c_{i_2}(2) 
   \cdots c_{i_n}(n)
\]
and the double universal Schubert polynomial by
\[ \Schub_w(c;d) = \sum_{u\cdot v = w} (-1)^{\ell(u)}
   \Schub_{u^{-1}}(d) \Schub_v(c) \,. 
\]
Since the ordinary Schubert polynomial $\Schub_w(X)$ does not depend
on which symmetric group $w$ belongs to, the same is true for
$\Schub_w(c;d)$.

As explained in the introduction, the double universal Schubert
polynomials describe the degeneracy loci $\Omega_w(G_\bull \to
F_\bull)$ of morphisms between vector bundles on a smooth algebraic
variety $\X$. The precise statement is the following result.

\begin{thm}[Fulton \cite{F2}] \label{thm:fulton}
  If the codimension of $\Omega_w(G_\bull \to F_\bull)$ is equal to
  $\ell(w)$ (or if this locus is empty) then the class of
  $\Omega_w(G_\bull \to F_\bull)$ in the cohomology ring of $\X$ is
  obtained from $\Schub_w(c;d)$ by evaluating at the Chern classes of
  the bundles, i.e., setting $c_i(j) = c_i(F_j)$ and $d_i(j) =
  c_i(G_j)$.
\end{thm}

We will need the following consequence of Theorem \ref{thm:fulton},
which generalizes a result of Kirillov~\cite[(5.5)]{K}, see
also~\cite[(26)]{F2}.

\begin{prop} \label{prop:same}
  Choose $m\gequ 0$ and substitute $b(j)$ for $c(j)$ and for $d(j)$ in
  $\Schub_w(c;d)$ for all $j \gequ m+1$.  We then have
\begin{multline*}
  \Schub_w(c(1),\dots,c(m),b(m+1),\dots \,;\,
  d(1),\dots,d(m),b(m+1),\dots) \\
  = \begin{cases} \Schub_w(c;d) & \text{if $w \in S_{m+1}$} \\
  0 & \text{otherwise.} \end{cases}
\end{multline*}
\end{prop}
\begin{proof}
  If $w \in S_{m+1}$ then it follows from the definition that
  $\Schub_w(c;d)$ is independent of $c_i(j)$ and $d_i(j)$ for $j \gequ
  m+1$.
  
  Assume that $w \in S_{n+1} \smallsetminus S_n$ where $n > m$.  We
  claim that $\Schub_w(c;d)$ vanishes as soon as we set $c_i(n) =
  d_i(n)$.  To see this, consider a variety $\X$ with bundles $F_j$
  for $1 \lequ j \lequ n$ and $G_j$ for $1 \lequ j \lequ n-1$, such
  that $\rank(F_j) = \rank(G_j) = j$, and such that all monomials of
  degree $\ell(w)$ in the Chern classes of these bundles are linearly
  independent.  For example we can take a product of Grassmann
  varieties $\X = \Gr(1,N) \times \dots \times \Gr(n,N) \times \dots
  \times \Gr(1,N)$ for $N$ large and use the tautological bundles.
  
  If we set $G_n = F_n$ then we have a sequence of bundles
  (\ref{eqn:univseq}) for which the map $G_n \to F_n$ is the identity
  and all other maps are zero.  Since $r_w(n,n) = n-1$ it follows that
  the locus $\Omega_w(G_\bull \to F_\bull)$ is empty.
  \refthm{thm:fulton} therefore implies that $\Schub_w(c;d)$ is zero
  when $c_i(n) = d_i(n) = c_i(F_n)$ for $1 \lequ i \lequ n$.
\end{proof}

The next identity is due to Kirillov \cite{K}; we 
include a proof for completeness.

\begin{cor}[Kirillov] \label{cor:cauchy}
  If $b_i(j)$, $c_i(j)$, and $d_i(j)$ are three sets of variables,
  \linebreak $1\lequ i \lequ j \lequ n$, then we have
\[ \Schub_w(c;d) = \sum_{u \cdot v = w} \Schub_u(b;d) \Schub_v(c;b) \,.\]
\end{cor}
\begin{proof}
  Let $\Schub(c;d)$ denote the function from permutations to
  polynomials which maps $w$ to $\Schub_w(c;d)$.  Following \cite[\S
  6]{M1} we define the product of two such functions $f$
  and $g$ by the formula
\[ (fg)(w) = \sum_{u\cdot v=w} f(u)g(v) \,. \]
Fulton's definition of universal Schubert polynomials then says that
$\Schub(c;d) = \Schub(0;d) \Schub(c;0)$, and \refprop{prop:same} shows
that $\Schub(0;c) \Schub(c;0) = \Schub(c;c) = \uu{1}$, where
$\uu{1}(w) = \delta_{1,w}$.  Now (6.6) of \cite{M1}
implies that also $\Schub(c;0)\Schub(0;c) = \uu{1}$.  We therefore
obtain 
$$\Schub(c;d) = \Schub(0;d) \Schub(c;0) = \Schub(0;d)
\Schub(b;0) \Schub(0;b) \Schub(c;0) = \Schub(b;d) \Schub(c;b),$$ 
as required.
\end{proof}

Later, we will need the following observation: \refprop{prop:same} and
\refcor{cor:cauchy} together imply that
\begin{multline} 
\label{eqn:splitusp1}
\Schub_w(c;d) = \sum_{u\cdot v = w}
  \Schub_u(0,\dots,0,c(r+1),c(r+2),\dots \,;\, d) \cdot
  \Schub_v(c(1),\dots,c(r);0)
\end{multline}
and similarly,
\begin{multline} \label{eqn:splitusp2}
\Schub_w(c;d) = \sum_{u\cdot v = w}
  \Schub_u(0;d(1),\dots,d(r)) \cdot
  \Schub_v(c; 0,\dots,0,d(r+1),d(r+2),\dots) \,.
\end{multline}

\subsection{Symmetric functions}
For each integer partition $\al=(\al_1\gequ\cdots\gequ\al_p\gequ 0)$, 
let $|\al|=\sum\al_i$ and let $\al'$ denote the conjugate of $\al$.  
Let $c=\{c_1,c_2,\ldots\}$ and $d=\{d_1,d_2,\ldots\}$ be ordered sets 
of independent variables.  Define the Schur determinant
\[
s_{\al}(c-d) = \det(h_{\al_i+j-i})_{p\times p}
\]
where the elements $h_k$ are determined by the identity of formal
power series
\[
\sum_{k\in \Z} h_k t^{k}=\frac{1-d_{1} t +d_{2} t^2 -\cdots}{1-c_{1} t +
c_2 t^2 -\cdots} \,. 
\]
In particular, $h_{0}=1$ and $h_{k}=0$ for $k<0$. The {\em
  supersymmetric Schur functions} $s_{\al}(X/Y)$ are obtained by
setting $c_i=e_i(X)$ and $d_i=e_i(Y)$ for all $i$, and the usual Schur
polynomials are given by the specializations
\[
s_{\al}(X/Y)\vert_{Y=0}=s_{\al}(X) \qquad \mathrm{and}
\qquad
s_{\al}(X/Y)\vert_{X=0}=(-1)^{|\al|}s_{\al'}(Y).
\]
If $E$ and $F$ are two vector bundles with total Chern classes $c(E)$
and $c(F)$, respectively, we will denote $s_{\al}(c(E)-c(F))$ by
$s_{\al}(E-F)$.

Let $\Lambda$ denote the ring of symmetric functions (as in
\cite{M3}).  For each permutation $w \in S_n$ there is a stable
Schubert polynomial or Stanley symmetric function $F_w \in \Lambda$
which is uniquely determined by the property that
\[
F_w(x_1,\ldots,x_k)=\AS_{1^m\times w}(x_1,\ldots,x_k)
\]
for all $m\gequ k$. \footnote{In Stanley's notation, the function
$F_{w^{-1}}$ is assigned to $w$.}  Here $1^m \times w \in S_{n+m}$ is
the permutation which is the identity on $\{1,\dots,m\}$ and which
maps $j$ to $w(j-m)+m$ for $j > m$ (see \cite[(7.18)]{M1}).  When
$F_w$ is written in the basis of Schur functions, one has
\[
F_w=\sum_{\al\, : \, |\al|=\ell(w)} d_{w\al} s_{\al}
\]
for some nonnegative integers $d_{w\al}$ \cite{EG,LS2}.  Fomin and
Greene show that the coefficient $d_{w\al}$ equals the number of
semistandard tableaux $T$ of shape $\al'$ such that the column word of
$T$ is a reduced word for $w$ \cite[Thm.~1.2]{FG}.

\section{Quiver varieties}
\label{qv}

\subsection{Definitions}
Let $E_1 \to E_2 \to \dots \to E_n$ be a sequence of vector bundles
and bundle maps over a non-singular variety $\X$.  Given rank conditions
$r = \{ r_{ij} \}$ for $1 \lequ i < j \lequ n$ there is a {\em quiver
  variety\/} given by
\[ \Omega_r(E_\bull) = \{ x \in \X \mid 
   \rank(E_i(x) \to E_j(x)) \lequ r_{ij} ~\forall i<j \} \,.
\]
For convenience, we set $r_{ii} = \rank E_i$ for all $i$, and we
demand that the rank conditions satisfy $r_{ij} \gequ \max(r_{i-1,j},
r_{i,j+1})$ and $r_{ij} + r_{i-1,j+1} \gequ r_{i-1,j} + r_{i,j+1}$ for
all $i \leq j$.
In this case, the expected codimension of $\Omega_r(E_\bull)$ is the
number $d(r) = \sum_{i<j} (r_{i,j-1}-r_{ij})(r_{i+1,j}-r_{ij})$.  The
main result of \cite{BF} states that when the quiver variety
$\Omega_r(E_\bull)$ has this codimension, its cohomology class is
given by
\begin{equation} 
\label{eqn:quiverformula}
  [\Omega_r(E_\bull)] = \sum_\lambda c_\lambda(r)\, 
  s_{\lambda^1}(E_2-E_1) \cdots s_{\lambda^{n-1}}(E_n-E_{n-1}) \,.
\end{equation}
Here the sum is over all sequences of partitions $\lambda =
(\lambda^1, \dots, \lambda^{n-1})$ such that $\sum |\lambda^i| =
d(r)$, and the coefficients $c_\lambda(r)$ are integers computed by a
combinatorial algorithm which we will not reproduce here.
\footnote{Knutson, Miller, and Shimozono have recently announced that 
they can prove that the integers $c_{\lambda}(r)$ are nonnegative.} 
These coefficients are uniquely determined by the condition that
(\ref{eqn:quiverformula}) is true for all varieties $\X$ and sequences
$E_\bull$, as well as the condition that $c_\lambda(r) =
c_\lambda(r')$, where $r'=\{r'_{ij}\}$ is the set of rank conditions
given by $r'_{ij} = r_{ij}+1$ for all $i \lequ j$.

The loci associated with universal Schubert polynomials are special
cases of these quiver varieties.  Given $w \in S_{n+1}$ we define rank
conditions $r^{(n)} = \{ r^{(n)}_{ij} \}$ for $1 \lequ i \lequ j \lequ
2n$ by
\[ r^{(n)}_{ij} = \begin{cases}
   r_w(2n+1-j, i) & \text{if $i\lequ n < j$} \\
   i & \text{if $j \lequ n$} \\
   2n+1-j & \text{if $i \gequ n+1$.}
\end{cases} \]
Then $\Omega_w(G_\bull \to F_\bull)$ is identical to the quiver
variety $\Omega_{r^{(n)}}(G_\bull \to F_\bull)$, and furthermore we
have $\ell(w) = d(r)$.  If we let $c^{(n)}_{w,\lambda} =
c_\lambda(r^{(n)})$ denote the quiver coefficients corresponding to
this locus, it follows that
\begin{equation} 
\label{eqn:quiv2univ}
 \Schub_w(c;d) = \sum_\lambda c_{w,\lambda}^{(n)}\, 
  s_{\lambda^1}(d(2)-d(1)) \cdots
  s_{\lambda^n}(c(n)-d(n)) \cdots s_{\lambda^{2n-1}}(c(1)-c(2)) \,.
\end{equation}

\subsection{Proof of Theorem \ref{introthm}}
It will be convenient to work with the element $\linebreak$
$P^{(n)}_w \in \Lambda^{\otimes 2n-1}$ defined by
\[ P^{(n)}_w = \sum_{\lambda} c_{w,\lambda}^{(n)}\, 
   s_{\lambda^1} \otimes \dots \otimes s_{\lambda^{2n-1}} \,.
\]
Theorem \ref{introthm} is a consequence of Fomin and Greene's formula
for stable Schubert polynomials combined with the following result.

\begin{thm} \label{thm:stanprd}
  For $w \in S_{n+1}$ we have
  \[ P^{(n)}_w = \sum_{u_1 \dots u_{2n-1} = w} 
     F_{u_1} \otimes \dots \otimes F_{u_{2n-1}}
  \]
  where the sum is over all reduced factorizations $w = u_1 \cdots
  u_{2n-1}$ such that $u_i \in S_{\min(i,2n-i)+1}$ for each $i$.
\end{thm}
\begin{proof}
  Since $r_w(p,q)+m = r_{1^m\times w}(p+m,q+m)$ for $m \geq 0$, it
  follows that the coefficients $c^{(n)}_{w,\lambda}$ are uniquely
  defined by the condition that
\begin{multline} \label{eqn:shift2quiv}
  \Schub_{1^m\times w}(c;d) = \sum_\lambda c_{w,\lambda}^{(n)}\,
  s_{\lambda^1}(d(2+m)-d(1+m)) \cdots \\ s_{\lambda^n}(c(n+m)-d(n+m)) 
  \cdots s_{\lambda^{2n-1}}(c(1+m)-c(2+m))
\end{multline}
for all $m \gequ 0$ (see also \cite[\S 4]{B}).

Given any two integers $p \lequ q$ we let $P_w^{(n)}[p,q]$ denote
the sum of the terms of $P^{(n)}_w$ for which $\lambda^i$ is empty when
$i<p$ or $i>q$:
\begin{equation*}
  P^{(n)}_w[p,q] = \sum_{\lambda: \lambda^i = \emptyset 
  \text{ for } i \not \in [p,q]} 
  c_{w,\lambda}^{(n)}\, s_{\lambda^1} \otimes \dots \otimes 
  s_{\lambda^{2n-1}} \,.
\end{equation*}

\begin{lemma} \label{lemma:splitquiv}
For any $1 < i \lequ 2n-1$ we have
\begin{equation}
\label{lem1eq} 
P^{(n)}_w = 
   \sum_{u \cdot v = w} P^{(n)}_u[1,i-1]\cdot P^{(n)}_v[i,2n-1] \,. 
\end{equation}
\end{lemma}
\begin{proof}
  We will do the case $i \lequ n$; the other one is similar.  For $f =
  \sum c_\lambda s_{\lambda^1} \otimes \dots \otimes
  s_{\lambda^{2N-1}} \in \Lambda^{\otimes 2N-1}$, we set
\[ f(c;d) = \sum c_\lambda\, s_{\lambda^1}(d(2)-d(1)) \cdots
   s_{\lambda^N}(c(N)-d(N)) \cdots s_{\lambda^{2N-1}}(c(1)-c(2)) \,.
\]
Equation (\ref{eqn:shift2quiv}) implies that $P^{(n)}_w \in
\Lambda^{\otimes 2n-1}$ is the unique element satisfying that
$(1^{\otimes m} \otimes P^{(n)}_w \otimes 1^{\otimes m}) (c;d) =
\Schub_{1^m\times w}(c;d)$ for all $m$.  This uniqueness is preserved
even if we set $d(i+m) = 0$.  The right hand side of the identity 
(\ref{lem1eq}) satisfies this by equation (\ref{eqn:splitusp2}) 
applied to $1^m \times w$.
\end{proof}

\begin{lemma} \label{lemma:stable}
For $1 \lequ i \lequ 2n-1$ we have
\[ P_w^{(n)}[i,i] = \begin{cases}
   1^{\otimes i-1} \otimes F_w \otimes 1^{\otimes 2n-1-i}
   & \text{if $w \in S_{m+1}$, $m=\min(i,2n-i)$} \\
   0 & \text{otherwise.}
\end{cases} \]
\end{lemma}
\begin{proof}
  If $w \not \in S_{m+1}$ then this follows from \refprop{prop:same},
  so assume $w \in S_{m+1}$.  For simplicity we will furthermore
  assume that $m = i$.  It is proved in \cite[\S 4]{B} that
  $P_w^{(m)}[m,m] = 1^{\otimes m-1} \otimes F_w \otimes 1^{\otimes
    m-1}$.  Let $\Phi : \Lambda^{\otimes 2m-1} \to \Lambda^{\otimes
    2n-1}$ be the linear map given by
\[ \Phi(s_{\lambda^1} \otimes \dots \otimes s_{\lambda^{2m-1}}) =
   s_{\lambda^1} \otimes \dots \otimes s_{\lambda^{m-1}} \otimes
   \Delta^{2n-2m}(s_{\lambda^m}) \otimes
   s_{\lambda^{m+1}} \otimes \dots \otimes s_{\lambda^{2m-1}}
\]
where $\Delta^{2n-2m} : \Lambda \to \Lambda^{\otimes 2n-2m+1}$ denotes
the $2n-2m$ fold coproduct, that is
\[ \Delta^{2n-2m}(s_{\lambda^m}) =
   \sum_{\tau_1,\dots,\tau_{2n-2m+1}}
   c^{\lambda^m}_{\tau_1,\dots,\tau_{2n-2m+1}}\,
   s_{\tau_1} \otimes \dots \otimes s_{\tau_{2n-2m+1}}
\]
where $c^{\lambda^m}_{\tau_1,\dots,\tau_{2n-2m+1}}$ is the coefficient
of $s_{\lambda^m}$ in the product $s_{\tau_1} s_{\tau_2} \cdots
s_{\tau_{2n-2m+1}}$.  For the definition of the locus
$\Omega_w(G_\bull \to F_\bull)$, the bundles $F_i$ and $G_i$ for $i
\gequ m+1$ are inessential in the sense of \cite[\S 4]{BF}.  Equation
(4.2) of loc.\ cit.\ therefore implies that $\Phi(P^{(m)}_w) =
P^{(n)}_w$.  Now the result follows from the identity
$P^{(n)}_w[m,2n-m] = \Phi(P^{(m)}_w[m,m]) = 1^{\otimes m-1} \otimes
\Delta^{2n-2m}(F_w) \otimes 1^{\otimes m-1}$.
\end{proof}

Theorem \ref{thm:stanprd} follows immediately from lemmas
\ref{lemma:splitquiv} and \ref{lemma:stable}.
\end{proof}

\begin{exa}
  For the permutation $w=s_{2}s_{1}=312$ in $S_3$,
the sequences of tableaux which satisfy the conditions of 
Theorem~\ref{introthm} are 
\begin{equation}\nonumber
(\ \emptyset \ ,\
\begin{picture}(10,20)
\put(0,15){\line(1,0){10}}
\put(0,15){\line(0,-1){20}}
\put(0,5){\line(1,0){10}}
\put(10,-5){\line(0,1){10}}
\put(0,-5){\line(1,0){10}}
\put(10,5){\line(0,1){10}}
\put(3,7){$1$}
\put(3,-3){$2$}
\end{picture}
\ ,\ \emptyset \ ) 
\ \ \mbox{ and } \ \ 
(\ \emptyset \ ,\ 
\begin{picture}(10,10)
\put(0,-2){\line(1,0){10}}
\put(0,-2){\line(0,1){10}}
\put(10,-2){\line(0,1){10}}
\put(0,8){\line(1,0){10}}
\put(3,0){$2$}
\end{picture} \ , \
\begin{picture}(10,10)
\put(0,-2){\line(1,0){10}}
\put(0,-2){\line(0,1){10}}
\put(10,-2){\line(0,1){10}}
\put(0,8){\line(1,0){10}}
\put(3,0){$1$}
\end{picture} \ ) \,.
\end{equation}
It follows that
\begin{eqnarray*}
\AS_{312}(c;d) & = & s_2(c(2)-d(2)) + s_{1}(c(2)-d(2))s_1(c(1)-c(2)) \\ 
& = & c_{1}(1)c_{1}(2)-c_{1}(1)d_{1}(2)-c_{2}(2) + d_2(2) \,.
\end{eqnarray*}
\end{exa}

In \cite{BF}, a conjectural combinatorial rule for general quiver
coefficients $c_{\lambda}(r)$ was given.  Although this rule was also
stated in terms of sequences of semistandard tableaux satisfying
certain conditions, it is different from Theorem \ref{introthm} in the
case of universal Schubert polynomials.  It would be interesting to
find a bijection that establishes the equivalence of these two rules.

\subsection{Skipping bundles} 
\label{skbs}

A permutation $w$ has a descent position at $i$ if $w(i)>w(i+1)$.
We say that a sequence $\{a_k\}\, :\, a_1<\cdots<a_p$ of integers is
{\em compatible with} $w$ if all descent positions of $w$ are
contained in $\{a_k\}$. 
Suppose that $w \in S_{n+1}$ and let $1 \lequ a_1 < a_2 < \dots < a_p
\lequ n$ and $1 \lequ b_1 < b_2 < \dots < b_q \lequ n$
be two sequences compatible with $w$ 
and $w^{-1}$, respectively. 

We let $E_\bull$ denote the subsequence
\[G_{b_1} \to G_{b_2} \to \dots \to G_{b_q} \to
F_{a_p} \to \dots \to F_{a_2} \to F_{a_1},\] 
and define rank conditions
$\tilde r^{(n)} = \{ \tilde r^{(n)}_{ij} \}$ for $1 \lequ i \lequ j \lequ
p+q$ by
\[ \tilde r^{(n)}_{ij} = \begin{cases}
  r_w(a_{p+q+1-j}, b_i) & \text{if $i \lequ q < j$} \\
  b_i & \text{if $j \lequ q$} \\
  a_{p+q+1-j} & \text{if $i \gequ q+1$}
\end{cases} \]
  Then the expected codimension of the locus $\Omega_{\tilde
  r^{(n)}}(E_\bull)$ is equal to $\ell(w)$.  However, in general this
  locus may contain $\Omega_w(G_\bull \to F_\bull)$ as a proper closed
  subset.  We will need the following criterion for equality (see also
  the remarks in \cite[\S 3]{F2} and Exercise 10 of \cite[\S 10]{F3}.)

\begin{lemma} \label{lemma:skip}
  Suppose that the map $G_{i-1} \to G_i$ is injective for $i \not \in
  \{b_k\}$ and the map $F_i \to F_{i-1}$ is surjective for $i \not \in
  \{a_k\}$.  Then $\Omega_{\tilde r^{(n)}}(E_\bull) = \Omega_w(G_\bull
  \to F_\bull)$ as subschemes of $\X$.
\end{lemma}
\begin{proof}
  Let $1 \lequ i,j \lequ n$ be given.  If $i$ is not a descent position
  for $w$ then either $w(i) \lequ j$ or $w(i+1) > j$.  In the first
  case this implies that $r_w(i,j) = r_w(i-1,j)+1$ so the condition
  $\rank(G_j \to F_i) \lequ r_w(i,j)$ follows from $\rank(G_j \to
  F_{i-1}) \lequ r_w(i-1,j)$ because the map $F_i \to F_{i-1}$ is
  surjective.  In the second case we have $r_w(i,j) = r_w(i+1,j)$ so
  the rank condition on $G_j \to F_i$ follows from the one on $G_j \to
  F_{i+1}$.  A similar argument works if $w^{-1}$ does not have a
  descent at position $j$.  We conclude that the locus
  $\Omega_w(G_\bull \to F_\bull)$ does not change if the bundles $G_j$
  for $j \not \in \{b_k\}$ and $F_i$ for $i \not \in \{a_k\}$ are
  disregarded.
\end{proof}

\begin{cor}
\label{relcor}
Let $w \in S_{n+1}$ and $\{a_k\}$ and $\{b_k\}$ be as above.  Then we
have
\[ [\Omega_{\tilde r^{(n)}}(E_\bull)] = \sum_\mu 
   \tilde c^{(n)}_{w,\mu}\, s_{\mu^1}(G_{b_2} - G_{b_1}) \cdots
   s_{\mu^q}(F_{a_p} - G_{b_q}) \cdots s_{\mu^{p+q-1}}(F_{a_1} -
   F_{a_2})
\]
with coefficients $\tilde c_{w,\mu}^{(n)} = c^{(n)}_{w,\lambda}$ where
the sequence $\lambda = (\lambda^1,\ldots,\lambda^{2n-1})$ is given by
\[ \lambda^i = \begin{cases} \mu^k & \text{if $i = b_k$} \\
  \mu^{p+q-k} & \text{if $i = 2n-a_k$} \\
  \emptyset & \text{otherwise.}
  \end{cases}
\]
\end{cor}
\begin{proof}
  We may assume that $F_i = F_{i-1} \oplus \C$ when $i \not \in
  \{a_k\}$ and $G_j = G_{j-1} \oplus \C$ when $j \not \in \{b_k\}$.
  In this case, notice that $s_\al(F_i - F_{i-1})$ is non-zero only
  when $\al$ is the empty partition or when $i = a_k$ for some $k$.
  Similarly $s_\al(G_{i-1}-G_i)$ is zero unless $\al$ is empty or $i =
  b_k$ for some $k$.  The result therefore follows from
  \reflemma{lemma:skip} and equation (\ref{eqn:quiv2univ}).
\end{proof}

\section{Schubert polynomials}
\label{sp}

\subsection{Degeneracy loci}
In this section, we will interpret the previous results for ordinary
double Schubert polynomials.  Let $V$ be a vector bundle of rank $n$
and let
\[ G_1 \subset G_2 \subset \dots \subset G_{n-1} \subset V 
   \tto F_{n-1} \tto \cdots \tto F_2 \tto F_1
\]
be a complete flag followed by a dual complete flag of $V$.  If $w \in
S_n$ then Fulton has proved \cite{F1} that
\[ [\Omega_w(G_\bull \to F_\bull)] =
   \Schub_w(x_1,\dots,x_n; y_1,\dots,y_n)
\]
where $x_i = c_1(\ker(F_i \to F_{i-1}))$, $y_i = c_1(G_i/G_{i-1})$,
and $\Schub_w(X;Y)$ is the double Schubert polynomial of Lascoux and
Sch{\"u}tzenberger.


Set $G'_i = V/G_i$ and $F'_i = \ker(V \to F_i)$.  Then we have a
sequence
\[ F'_{n-1} \subset \dots \subset F'_1 \subset V 
   \tto G'_1 \tto \cdots \tto G'_{n-1}
\]
and it is easy to check that $\Omega_w(G_\bull \to F_\bull) =
\Omega_{w_0 w^{-1} w_0}(F'_\bull \to G'_\bull)$ as subschemes of $\X$,
where $w_0$ is the longest permutation in $S_n$.

Let $1 \lequ a_1 < \dots < a_p \lequ n-1$ and 
$0 \lequ b_1 < \dots <b_q \lequ n-1$ be two sequences compatible with $w$ 
and $w^{-1}$, respectively.
Then by applying section \ref{skbs} to the
subsequence $F'_{a_p} \to \dots \to F'_{a_1} \to G'_{b_1} \to \dots \to
G'_{b_q}$ we obtain
\begin{multline*} 
  [\Omega_w(G_\bull \to F_\bull)] = \\
  \sum_\mu \tilde c_{w_0 w^{-1} w_0,\mu}^{(n-1)}\,
  s_{\mu^1}(F'_{a_{p-1}} - F'_{a_p}) \cdots s_{\mu^p}(G'_{b_1}
  - F'_{a_1}) \cdots s_{\mu^{p+q-1}}(G'_{b_q} - G'_{b_{q-1}}) \,.
\end{multline*}

Set $a_0 = b_0 = 0$.  If we let $X_i = \{x_{a_{i-1}+1}, \dots,
x_{a_i}\}$ denote the Chern roots of $\ker(F_{a_i} \to F_{a_{i-1}})$
and $Y_i = \{y_{b_{i-1}+1}, \dots, y_{b_i}\}$ be the Chern roots of
$G_{b_i}/G_{b_{i-1}}$ then the previous equality can be written as
\begin{equation} \label{eqn:ss}
  \Schub_w(X;Y) = \sum_\mu \tilde c_{w_0 w^{-1} w_0,\mu}^{(n-1)}\,
  s_{\mu^1}(X_p) \cdots s_{\mu^p}(X_1/Y_1) \cdots
  s_{\mu^{p+q-1}}(0/Y_q) \,.
\end{equation}

This equation is true in the cohomology ring $H^*(\X;\Z)$, in which
there are relations between the variables $x_i$ and $y_i$ (including
e.g.\ the relations $e_j(x_1,\dots,x_n) = c_j(V) = e_j(y_1,\dots,y_n)$
for $1 \lequ j \lequ n$).  We claim, however, that (\ref{eqn:ss})
holds as an identity of polynomials in independent variables.  For
this, notice that the identity is independent of $n$, i.e.\ the
coefficient $\tilde c^{(n-1)}_{w_0 w^{-1} w_0,\mu}$ does not change when $n$
is replaced with $n+1$ and $w_0$ with the longest element in
$S_{n+1}$.  If we choose $n$ sufficiently large, we can construct a
variety $\X$ on which (\ref{eqn:ss}) is true, and where all monomials
in the variables $x_i$ and $y_i$ of total degree at most $\ell(w)$ are
linearly independent, which establishes the claim.

\subsection{Splitting Schubert polynomials}

We continue by reformulating equation (\ref{eqn:ss}) to
obtain a more natural expression for double Schubert polynomials.

It follows from (\ref{eqn:ss}) together with the main result of
\cite{B} that $F_{w_0 w^{-1} w_0} = F_w$.  Therefore, the coefficient
$d_{w \alpha}$ of the Schur expansion of $F_w$ is equal to the number
of semistandard tableaux of shape $\alpha'$ such that the column word
is a reduced word for $w_0 w^{-1} w_0$.  Notice that if $e = (e_1,
e_2, \dots, e_\ell)$ is a reduced word for $w_0 w^{-1} w_0$ then $\wt
e = (n+1-e_\ell, \dots, n+1-e_1)$ is a reduced word for $w$.
Furthermore, if $e$ is the column word of a tableau of shape
$\alpha'$, then $\wt e$ is the column word of a skew tableau whose
shape is the 180 degree rotation of $\alpha'$.  If we denote this
rotated shape by $\widetilde \alpha$ then we conclude that $d_{w
  \alpha}$ is also equal to the number of skew tableaux of shape $\wt
\alpha$ such that the column word is a reduced word for $w$.  Theorem
\ref{thm:stanprd} now implies the following variation of our main
theorem:

\begin{thm1a}
  If $w \in S_{n+1}$ then the coefficient $c_{w,\lambda}^{(n)}$ equals
  the number of sequences of semistandard
  skew tableaux $(T_{1},\ldots,T_{2n-1})$
  such that $T_i$ has shape $\wt \lambda^i$, the entries of $T_{i}$
  are at most ${\rm min}(i,2n-i)$, and
  $\col(T_{1})\cdots\col(T_{2n-1})$ is a reduced word for $w$.
\end{thm1a}

We say that a sequence of tableaux $(T_1,\dots,T_r)$ is {\em
strictly bounded below\/} by an integer sequence $(a_1,\dots,a_r)$ if
the entries of $T_i$ are strictly greater than $a_i$, for each~$i$.

\begin{thm} \label{thm:splitnice}
  Let $w \in S_n$ and let $1 \lequ a_1 < \dots < a_p$ and $0 \lequ b_1 <
  \dots < b_q$ be two sequences compatible with
  $w$ and $w^{-1}$, respectively.  Then we have
\begin{equation} \label{eqn:splitnice}
   \Schub_w(X;Y) = \sum_\lambda c_\lambda \,
   s_{\lambda^1}(0/Y_q) \cdots s_{\lambda^q}(X_1/Y_1) \cdots 
   s_{\lambda^{p+q-1}}(X_p)
\end{equation}
where $X_i = \{x_{a_{i-1}+1},\dots,x_{a_i}\}$ and $Y_i =
\{y_{b_{i-1}+1},\dots,y_{b_i}\}$ and the sum is over all sequences of
partitions $\lambda = (\lambda^1,\dots,\lambda^{p+q-1})$.  Each
$c_\lambda$ is a quiver coefficient, equal to the number of sequences
of semistandard tableaux $(T_1,\dots,T_{p+q-1})$ strictly bounded
below by $(b_{q-1}, \dots, b_1, 0, a_1, a_2, \dots, a_{p-1})$, such
that the shape of $T_i$ is conjugate to $\lambda^i$ and
$\col(T_1)\cdots \col(T_{p+q-1})$ is a reduced word for $w$.
\end{thm}
\begin{proof}
  This follows from equation (\ref{eqn:ss}) together with Theorem $1'$
  applied to \linebreak $w_0 w^{-1} w_0$.  To translate between the
  sequences of skew tableaux in Theorem $1'$ and the sequences in the
  present theorem, simply rotate a whole sequence of skew tableaux by
  180 degrees (this means invert the order of the sequence, and turn
  each skew tableau on its head).  Then replace each entry $e$ with
  $n+1-e$.
  
  A purely algebraic proof of the splitting formula (\ref{eqn:splitnice})
  is also possible. This requires versions of the results of Section
  \ref{usp} for ordinary Schubert polynomials, available e.g.\ from
  \cite{L2} and \cite{FK}, together with \cite[Thm.\ 1.2]{FG}.
  However, some geometric reasoning is needed to interpret the
  $c_{\la}$ as quiver coefficients.
\end{proof}

Notice that if one takes $b_1 = 0$ in Theorem \ref{thm:splitnice} then
the set of variables $Y_1$ is empty, so equation (\ref{eqn:splitnice})
contains only signed products of single Schur polynomials.
Observe also that a factor $s_{\lambda^i}(0/Y_k)$ in
(\ref{eqn:splitnice}) will vanish if $\lambda^i$ has more than
$b_k-b_{k-1}$ columns and that $s_{\lambda^i}(X_k)$ vanishes if
$\lambda^i$ has more than $a_k-a_{k-1}$ rows.  Therefore equation
(\ref{eqn:splitnice}) uses only a subset of the quiver coefficients
of Theorem \ref{thm:splitnice}.

\begin{exa} Let $w = 321$ be the longest element in $S_3$, and 
choose the sequences $\{a_i\} = \{b_i\} = \{1 < 2\}$. Then the four
sequences of tableaux satisfying the 
conditions of Theorem~\ref{thm:splitnice} 
are 
\begin{equation}\nonumber
(\ \begin{picture}(10,10)
\put(0,-2){\line(1,0){10}}
\put(0,-2){\line(0,1){10}}
\put(10,-2){\line(0,1){10}}
\put(0,8){\line(1,0){10}}
\put(3,0){$2$}
\end{picture} \ , 
\ \begin{picture}(10,10)
\put(0,-2){\line(1,0){10}}
\put(0,-2){\line(0,1){10}}
\put(10,-2){\line(0,1){10}}
\put(0,8){\line(1,0){10}}
\put(3,0){$1$}
\end{picture} \ , 
\ \begin{picture}(10,10)
\put(0,-2){\line(1,0){10}}
\put(0,-2){\line(0,1){10}}
\put(10,-2){\line(0,1){10}}
\put(0,8){\line(1,0){10}}
\put(3,0){$2$}
\end{picture} \ )\,, \ \ \ 
(\ \begin{picture}(10,10)
\put(0,-2){\line(1,0){10}}
\put(0,-2){\line(0,1){10}}
\put(10,-2){\line(0,1){10}}
\put(0,8){\line(1,0){10}}
\put(3,0){$2$}
\end{picture} \ , \ 
\begin{picture}(20,10)
\put(0,-2){\line(1,0){20}}
\put(0,8){\line(1,0){20}}
\put(0,-2){\line(0,1){10}}
\put(10,-2){\line(0,1){10}}
\put(20,-2){\line(0,1){10}}
\put(3,0){$1$}
\put(13,0){$2$}
\end{picture} \ ,\  \emptyset \ )\,,\ \ 
(\ \emptyset \ , \
\begin{picture}(10,20)
\put(0,15){\line(1,0){10}}
\put(0,15){\line(0,-1){20}}
\put(0,5){\line(1,0){10}}
\put(10,-5){\line(0,1){10}}
\put(0,-5){\line(1,0){10}}
\put(10,5){\line(0,1){10}}
\put(3,7){$1$}
\put(3,-3){$2$}
\end{picture}
\ ,\ 
\begin{picture}(10,10)
\put(0,-2){\line(1,0){10}}
\put(0,-2){\line(0,1){10}}
\put(10,-2){\line(0,1){10}}
\put(0,8){\line(1,0){10}}
\put(3,0){$2$}
\end{picture} \ )\,, 
\ \ \mbox{ and } \ \ 
(\ \emptyset \ , \
\begin{picture}(20,20)
\put(0,15){\line(1,0){20}}
\put(0,15){\line(0,-1){20}}
\put(0,5){\line(1,0){20}}
\put(20,5){\line(0,1){10}}
\put(10,-5){\line(0,1){10}}
\put(0,-5){\line(1,0){10}}
\put(10,5){\line(0,1){10}}
\put(3,7){$1$}
\put(3,-3){$2$}
\put(13,7){$2$}
\end{picture} \ ,
\ \emptyset \ )\,,
\end{equation}
all of which give nonvanishing terms and correspond to the reduced
word $s_2 s_1 s_2$.  We thus have
\begin{gather*}
\AS_{321}(X;Y) = \\
s_1(0/y_2) s_1(x_1/y_1) s_1(x_2) + s_1(0/y_2)s_{1,1}(x_1/y_1) + 
s_2(x_1/y_1)s_1(x_2) + s_{2,1}(x_1/y_1) \\
= - y_2(x_1-y_1)x_2 - y_2(y_1^2-x_1y_1)
  + (x_1^2-x_1y_1)x_2 + (x_1y_1^2-x_1^2y_1) \\
= (x_1-y_1)(x_1-y_2)(x_2-y_1) \,.
\end{gather*}
\end{exa}

\begin{cor} 
\label{cor:singlesplitcor}
Suppose that $w\in S_n$ is a permutation compatible with the
sequence $a_1< \cdots < a_p$.  Then we have
\[
   \AS_w(X) = \sum_\lambda c_\lambda \,
   s_{\lambda^1}(X_1) \cdots s_{\lambda^p}(X_p)
\]
where $X_i = \{x_{a_{i-1}+1},\dots,x_{a_i}\}$ and
the sum is over all sequences of
partitions $\lambda = (\lambda^1,\dots,\lambda^p)$.  Each
$c_\lambda$ is a quiver coefficient, equal to the number of sequences
of semistandard tableaux $(T_1,\dots,T_p)$ strictly bounded
below by $(0, a_1, a_2, \dots, a_{p-1})$, such
that the shape of $T_i$ is conjugate to $\lambda^i$ and
$\col(T_1)\cdots \col(T_p)$ is a reduced word for $w$.
\end{cor}

\begin{exa} Consider the permutation $w=s_1s_2s_1s_3s_4s_3=32541$ in
  $S_5$, with descent positions at $1$, $3$, and $4$. The two sequences of
tableaux satisfying the conditions of
Corollary~\ref{cor:singlesplitcor} which give nonvanishing
terms are
\begin{equation}\nonumber
(\
\begin{picture}(10,30) 
\put(0,-12){\line(0,1){30}}
\put(0,-12){\line(1,0){10}}
\put(0,-2){\line(1,0){10}}
\put(0,8){\line(1,0){10}}
\put(0,18){\line(1,0){10}}
\put(10,-12){\line(0,1){30}}
\put(3,-10){$4$}
\put(3,0){$2$}
\put(3,10){$1$}
\end{picture} \ , \
\begin{picture}(20,10)
\put(0,-2){\line(1,0){20}}
\put(0,8){\line(1,0){20}}
\put(0,-2){\line(0,1){10}}
\put(10,-2){\line(0,1){10}}
\put(20,-2){\line(0,1){10}}
\put(3,0){$2$}
\put(13,0){$3$}
\end{picture} \ ,\
\begin{picture}(10,10)
\put(0,-2){\line(1,0){10}}
\put(0,-2){\line(0,1){10}}
\put(10,-2){\line(0,1){10}}
\put(0,8){\line(1,0){10}}
\put(3,0){$4$}
\end{picture} \ ) \ \ \mbox{ and } \ \  
(\ \begin{picture}(10,20)
\put(0,15){\line(1,0){10}}
\put(0,15){\line(0,-1){20}}
\put(0,5){\line(1,0){10}}
\put(10,-5){\line(0,1){10}}
\put(0,-5){\line(1,0){10}}
\put(10,5){\line(0,1){10}}
\put(3,7){$1$}
\put(3,-3){$2$}
\end{picture} \ , 
\ \begin{picture}(20,20)
\put(0,15){\line(1,0){20}}
\put(0,15){\line(0,-1){20}}
\put(0,5){\line(1,0){20}}
\put(20,5){\line(0,1){10}}
\put(10,-5){\line(0,1){10}}
\put(0,-5){\line(1,0){10}}
\put(10,5){\line(0,1){10}}
\put(3,7){$2$}
\put(3,-3){$4$}
\put(13,7){$3$}
\end{picture} \ , \
\begin{picture}(10,10)
\put(0,-2){\line(1,0){10}}
\put(0,-2){\line(0,1){10}}
\put(10,-2){\line(0,1){10}}
\put(0,8){\line(1,0){10}}
\put(3,0){$4$}
\end{picture} \ )\,.
\end{equation}

\medskip
\noin
The reduced words for $w$ corresponding to these sequences are 
$s_4s_2s_1s_2s_3s_4$ and $s_2s_1s_4s_2s_3s_4$, respectively. 
It follows that
\begin{eqnarray*}
\AS_{32541}(X)  &=&  s_3(x_1)s_{1,1}(x_2,x_3)s_1(x_4)+ 
s_2(x_1)s_{2,1}(x_2,x_3)s_1(x_4) \\
&=& x_1^3\cdot x_2x_3\cdot x_4 + x_1^2(x_2^2x_3+x_3^2x_2)x_4.
\end{eqnarray*}
\end{exa}

The special case of Corollary \ref{cor:singlesplitcor} with 
$a_k = k$ gives a
formula for the coefficient of each monomial in $\AS_w(X)$, which
is equivalent to that of \cite[Thm.\ 1.1]{BJS}. We deduce that these
monomial coefficients are quiver coefficients. The same conclusion
holds for double Schubert polynomials:

\begin{cor}
  Let $w \in S_n$ and let $x^u y^v = x_1^{u_1} \cdots
  x_{n-1}^{u_{n-1}} y_1^{v_1} \cdots y_{n-1}^{v_{n-1}}$ be a monomial
  of total degree $\ell(w)$.  Set $g_i = \sum_{k=n-i}^{n-1} v_k$ and
  $f_i = g_{n-1} + \sum_{k=1}^i u_k$.  Then the coefficient of $x^u
  y^v$ in the double Schubert polynomial $\Schub_w(X;Y)$ is equal to
  $(-1)^{g_{n-1}}$ times the number of reduced words
  $(e_1,\dots,e_{\ell(w)})$ for $w$ such that $n-i \lequ e_{g_{i-1}+1}
  < \dots < e_{g_i}$ and $e_{f_{i-1}+1} > \dots > e_{f_i} \gequ i$ for
  all $1 \lequ i \lequ n-1$.
\end{cor}
\begin{proof}
  Take $\{a_k\} = \{1,2,\dots,n-1\}$ and $\{b_k\} =
  \{0,1,2,\dots,n-1\}$ in Theorem \ref{thm:splitnice} and notice that
  a Schur polynomial $s_\alpha(0/y_i)$ is non-zero only if $\alpha =
  (1^{v_i})$ is a single column with $v_i$ boxes, in which case
  $s_\alpha(0/y_i) = (-y_i)^{v_i}$.  Similarly $s_\alpha(x_i)$ is
  equal to $x_i^{u_i}$ if $\alpha = (u_i)$ is a single row with $u_i$
  boxes, and is zero otherwise.  The reduced words of the corollary
  are exactly those that form a sequence of tableaux on the conjugates
  of these shapes and strictly bounded below by
  $(n-2,n-3,\dots,1,0,0,1,\dots,n-2)$.
\end{proof}

\begin{remark}
  Corollary~\ref{cor:singlesplitcor} implies that the single Schubert
  polynomial $\AS_w(X)$ is the character of a $\GL(a_1,\C) \times
  \GL(a_2-a_1,\C) \times \cdots \times \GL(a_p-a_{p-1},\C)$ module.
  Our expression for the coefficients $c_{\la}$ gives an isotypic
  decomposition of this module.  This is closely related to a theorem
  of Kraskiewicz and Pragacz \cite{KP}, which shows that $\AS_w(X)$ is
  a character of a Borel module.  P.~Magyar reports that the former
  property may be deduced from the latter (private communication).
\end{remark}

\section{Giambelli formulas}
\label{giamsec}

In this final section, we observe that Corollaries \ref{relcor} and
\ref{cor:singlesplitcor} give explicit, non-recursive answers to the
{\em Giambelli problem} for the classical and quantum cohomology of
partial flag manifolds.

Suppose that $V$ is a complex vector space of dimension $n$ and choose 
integers $0 = a_0 < a_1< \cdots < a_p<a_{p+1}=n$. Let $\X$ be the partial
flag variety which parametrizes quotients
\begin{equation}
\label{quots}
V\tto F_{a_p} \tto \cdots \tto F_{a_1}
\end{equation}
with $\rank(F_{a_k})=a_k$ for each $k$. We will also use (\ref{quots})
to denote the tautological sequence of quotient bundles over $\X$, and define
$Q_k=\mathrm{Ker}(F_{a_k}\to F_{a_{k-1}})$, for $1\lequ k\lequ p+1$.
According to Borel \cite{Bo}, the 
cohomology ring $H^*(\X;\Z)$ is presented as the
polynomial ring in the Chern classes $c_i(Q_k)$ for all $i$ and $k$, 
modulo the relation
\[ c(Q_1) c(Q_2) \cdots c(Q_{p+1}) = 1 \,. \]

Fix a complete flag $G_{\bull}$ of subspaces of $V$, and let $S(a)$ denote 
the subset of $S_n$ consisting of permutations $w$ compatible with $\{a_k\}$. 
For each $w\in S(a)$, there is a Schubert variety 
$\Omega_w\subset \X$, defined as the locus of $x\in\X$ such that 
\[
\rank(G_j(x)\ra F_i(x))\lequ r_w(i,j), \ \ \mathrm{for} \ \ 
 i\in \{a_1,\ldots a_p\} \ \ \mathrm{and}\ \ 1\lequ j\lequ n.
\]
The {\em Schubert classes} $[\Omega_w]\in H^{2\ell(w)}(\X;\Z)$ for
$w\in S(a)$ form a natural `geometric basis' for the cohomology ring of
$\X$.  The following Giambelli formula, which is a direct
consequence of Corollary \ref{cor:singlesplitcor}, writes these
classes as polynomials in the `algebraic generators' for $H^*(\X;\Z)$
given by the $c_i(Q_k)$.

\begin{giamI}
We have
\begin{equation}\label{giamI}
  [\Omega_w]= \sum_\lambda c_\lambda \,
  s_{\lambda^1}(Q_1) \cdots s_{\lambda^p}(Q_p), 
\end{equation}
where the sum is over sequences of partitions $\lambda =
(\lambda^1,\dots,\lambda^p)$ and $c_\lambda$ is the quiver coefficient
of Corollary \ref{cor:singlesplitcor}.
\end{giamI}

Alternatively, we can use Corollary \ref{relcor} to express the
Schubert class $[\Omega_w]$ as a polynomial in the special Schubert
classes $c_i(F_{a_k})$.  Notice that the two approaches are identical
for Grassmannians.


\begin{giamII}
We have
\begin{equation}\label{giamII}
[\Omega_w]= \sum_\nu 
  \widetilde{c}_{\nu}\, 
  s_{\nu^p}(F_{a_p}) \cdots s_{\nu^1}(F_{a_1} - F_{a_2}),
\end{equation}
where the sum is over sequences of partitions $\nu=(\nu^1,\ldots,\nu^p)$
and $\widetilde{c}_{\nu}=c^{(n)}_{w,\lambda}$, where $\lambda=(\lambda^1,
\ldots,\lambda^{2n-1})$ is given by
\[ \lambda^i = \begin{cases} 
  \nu^k & \text{if $i = 2n-a_{p+1-k}$} \\
  \emptyset & \text{otherwise.}
  \end{cases}
\]
\end{giamII}

When comparing the above two Giambelli formulas, recall that the the
quiver coefficients $c_{\lambda}$ which appear in (\ref{giamI})
correspond to the permutation $w_0w^{-1}w_0$.  The equivalence of
(\ref{giamI}) and (\ref{giamII}) can also be checked directly, by
using the Chern class identities $c(Q_k)=c(F_{a_k}-F_{a_{k-1}})$.

Following \cite{FGP} and \cite[\S 3.2]{C-F}, we recall that equation
(\ref{giamII}) can be used to obtain a quantum Giambelli formula which
holds in the small quantum cohomology ring $QH^*(\X)$.  For this, one
simply replaces all the special Schubert classes which appear in the
Schur determinants $s_{\nu^k}(F_{a_k} - F_{a_{k+1}})$ in
(\ref{giamII}) with the corresponding quantum classes, as in loc.\ 
cit.\ (compare also with \cite[Prop.\ 4.3]{F2}).

\end{document}